\documentclass{amsart}
\usepackage{graphicx}
\newtheorem{theorem}{Theorem}
\newtheorem{lemma}[theorem]{Lemma}
\newtheorem{corollary}[theorem]{Corollary}
\newtheorem{definition}{Definition}
\newtheorem{proposition}{Proposition}
 
\newcommand{\bla}{\bullet}
\newcommand{\whi}{\circ}

\newcommand{\basic}{\mathcal{B}}

\title{Crossing and alignments of permutations}

\author{Sylvie Corteel}
\address{CNRS LRI Universit\'e Paris-Sud B\^atiment 490 91405 Orsay Cedex}
\email{corteel@lri.fr}
\date{January 19th, 2006}
\begin{document}
\maketitle

\begin{abstract}
We introduce the notion of crossings and nestings of a permutation.
We compute the generating function of permutations with a fixed number of 
weak exceedances, crossings and nestings. 
We link alignments and permutation patterns to
these statistics.
We generalize to the case of decorated
permutations. Finally we show how this is related to
the stationary distribution of the Partially ASymmetric Exclusion Process (PASEP) model. 
\end{abstract} 
\section{Introduction}

We introduce the notion of crossings and nestings of a permutation.
The purpose of this paper is to link permutation patterns and
these statistics. Our main result
is the following:
\begin{theorem}
The number $B(n,k,\ell,m)$ of permutations $\sigma$ of 
$[n]$ with $k$ weak exceedances, $\ell$ crossings
and $m$ nestings is equal to  the
number $D(n,k,\ell,m)$ of permutations of $[n]$ with $n\mathrm{-}k$ descents, 
$\ell$ occurrences of the pattern $31\mathrm{-}2$ 
and $m$ occurrences of the pattern $2\mathrm{-}31$.
\label{max}
\end{theorem}
We prove this theorem  by exhibiting their generating function.
Let
\begin{equation}
F(q,p,x,y)=
\cfrac{1}{1-b_0x-\cfrac{\lambda_1x^2}{1-b_1x-\cfrac{\lambda_2x^2}
{1-b_2x-\cfrac{\lambda_3x^2}{\ddots}}}}
\label{best}
\end{equation}
with $[n]_{p,q}=p^{n-1}+p^{n-2}q+\cdots +pq^{n-2}+q^{n-1}$, $b_n=y[n+1]_{p,q}+[n]_{p,q}$ and $\lambda_n=y[n]_{p,q}^2$.
\begin{proposition}
The coefficient of $x^n y^k x^\ell y^m$ in $F(q,p,x,y)$ is equal to $B(n,k,\ell,m)$.
\label{th1}
\end{proposition}
\begin{proposition}
The coefficient of $x^n y^k x^\ell y^m$ in $F(q,p,x,y)$ is equal to  $D(n,k,\ell,m)$.
\label{th2}
\end{proposition}
We also propose a bijective proof of Theorem \ref{max}. To prove this we use 
a bijection of Foata and Zeilberger \cite{dv,fz},
a bijection due to Fran\c{c}on and Viennot \cite{fv}  and 
results from \cite{cm,csz}. \\

The notions of crossings and nestings are closely related
to the notion of alignments, defined in \cite{wi}, 
which  come from the enumeration
of totally positive Grassmann cells.  
The link between alignments and 
patterns was conjectured by Steingr\'imsson and Williams \cite{wi0}.

These alignments define some new $q$-analogs of the Eulerian numbers
$\hat{E}_{k,n}(q)$ were introduced 
by Williams \cite{wi} building on work of  Postnikov \cite{po}.
Let $[k]_q$ be $1+q+\cdots +q^{k-1}$, then 
\begin{proposition}\cite{wi}
The number of permutations of $[n]$ with $k$ weak exceedances
and $\ell$ alignments is the coefficient of
$q^{(k-1)(n\mathrm{-}k)-\ell}$ in 
$$
\hat{E}_{k,n}(q)=q^{k-k^2}\sum_{i=0}^{k-1}(-1)^i [k-i]_q^n q^{k(i-1)}\left( {n\choose i}q^{k-i}+{n\choose i-1}\right). 
$$
\end{proposition}
These numbers have the property that if $q=1$ they are the Eulerian numbers,
if $q=0$ they are the Narayana numbers and if $q=-1$ they are the binomial
coefficients. See \cite{wi} for details.
Williams proved in \cite{wi} using permutation tableaux
that ${E}(q,x,y)=\sum_{n,k}q^{n\mathrm{-}k}\hat{E}_{k,n}(q)y^k x^n$ is equal to
$$
\sum_{i=0}^{\infty} \frac{y^i (q^{2i+1} - y)}
{q^{i^2+i+1} (q^i - q^{i+1}[i]_qx + [i]_qxy)}.
$$                                                                                
Here we exhibit the generating function  
$\hat{E}(q,x,y)=\sum_{n,k}\hat{E}_{k,n}(q)y^k x^n$
and show that
\begin{theorem}
$$
\hat{E}(q,x,y)=F(q,1,x,y).
$$
\label{big}
\end{theorem}
This proves the conjecture of Steingr\'imsson and Williams \cite{wi0}.\\

We define the combinatorial statistics in Section 2.
We then prove Proposition \ref{th1}  and Theorem \ref{big}  in Section 3.
In Section 4 we prove Proposition \ref{th2}.  This also gives a proof
of  Theorem
\ref{max} which is implied by Propositions \ref{th1} and \ref{th2}.
In Section 5 we propose a direct bijective proof of Theorem
\ref{max}. 
In Section 6 we generalize to the case of decorated
permutations. In Section 7 we show how these numbers appear naturally
in the stationary distribution of the  PASEP model \cite{bcer}.

\section{Definitions}

Let $\sigma=(\sigma(1),\ldots ,\sigma(n))$ be a permutation
of $[n]=\{1,2,\ldots ,n\}$. The number of weak exceedances of a permutation
$\sigma$ is the cardinality of the set $\{j\ | \ \sigma(j)\ge j\}$. We denote
this number by $WEX(\sigma)$.

Let
\begin{itemize}
\item $C_+(i)=\{j\ | \ j<i\le\sigma(j)<\sigma(i)\}$,
\item $C_-(i)=\{j\ | \ j>i>\sigma(j)>\sigma(i)\}$, 
\end{itemize}
and $C_+(\sigma)=\sum_{i=1}^n |C_+(i)|$. 
\begin{definition}
The number of crossings of a permutation $\sigma$ is equal to
$$
C_+(\sigma)+C_-(\sigma).
$$
\end{definition}

\begin{figure}[ht!]
  \begin{center}
    \includegraphics[scale=0.5]{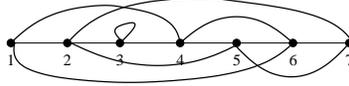}
  \end{center}\centering
  \caption{The permutation diagram of  $\sigma=(4,7,3,6,2,1,5)$}
  \label{perm}
\end{figure}

We can also define those parameters using the permutation diagram. 
 We draw a line and
put the numbers from 1 to $n$ and we draw an edge from
$i$ to $\sigma(i)$ above the line if $i\le \sigma(i)$ and 
under the line otherwise.
For example, the permutation diagram of  $\sigma=(4,7,3,6,2,1,5)$ is 
on Figure \ref{perm}. \\

$C_+(\sigma)$ is the number of pairs of
edges above the line that cross or touch \begin{picture}(0,0)%
\includegraphics{cplus.pstex}%
\end{picture}%
\setlength{\unitlength}{2072sp}%
\begingroup\makeatletter\ifx\SetFigFont\undefined%
\gdef\SetFigFont#1#2#3#4#5{%
  \reset@font\fontsize{#1}{#2pt}%
  \fontfamily{#3}\fontseries{#4}\fontshape{#5}%
  \selectfont}%
\fi\endgroup%
\begin{picture}(1460,380)(1973,-1464)
\end{picture}%
 and  $C_-(\sigma)$ is the number of pairs of
edges under the line that cross \begin{picture}(0,0)%
\includegraphics{cmoins.pstex}%
\end{picture}%
\setlength{\unitlength}{2072sp}%
\begingroup\makeatletter\ifx\SetFigFont\undefined%
\gdef\SetFigFont#1#2#3#4#5{%
  \reset@font\fontsize{#1}{#2pt}%
  \fontfamily{#3}\fontseries{#4}\fontshape{#5}%
  \selectfont}%
\fi\endgroup%
\begin{picture}(1460,380)(1973,-1018)
\end{picture}%
.

For $\sigma=(4,7,3,6,2,1,5)$ on Figure \ref{perm}, we have $C_{+}(\sigma)=2$ and $C_{-}(\sigma)=1.$
The pairs of edges
contributing to  $C_{+}(\sigma)$ are $\{(1,4),(4,6)\}$ and $\{(1,4),(2,7)\}$.
The pair of edges
contributing to  $C_{-}(\sigma)$ is $\{(7,5),(6,1)\}$.\\

Then let
\begin{itemize}
\item $A_+(i)=\{j\ |\ j<i\le \sigma(i)<\sigma(j)\}$;
\item $A_-(i)=\{j\ | \ j>i>\sigma(i)>\sigma(j)\}$;
\item $A_{+,-}(i)=
\{j\ |\ j\le \sigma(j)<\sigma(i)<i\}\cup \{j\ |\ \sigma(i)<i<j\le \sigma(j)\}$.
\end{itemize}
We then set 
$A_+(\sigma)=\sum_{i=1}^n |A_+(i)|$.

For $\sigma=(4,7,3,6,2,1,5)$, $A_+(\sigma)=3$, $A_-(\sigma)=1$, 
and $A_{+,-}(\sigma)=2$.\\

We can again define these parameters using the permutation diagram of
 permutation $\sigma$ of $[n]$. Then $A_+(\sigma)$ (resp. $A_-(\sigma)$) is the number of pairs of
nested edges above (resp. under) the line 
\begin{picture}(0,0)%
\includegraphics{aplus.pstex}%
\end{picture}%
\setlength{\unitlength}{2072sp}%
\begingroup\makeatletter\ifx\SetFigFont\undefined%
\gdef\SetFigFont#1#2#3#4#5{%
  \reset@font\fontsize{#1}{#2pt}%
  \fontfamily{#3}\fontseries{#4}\fontshape{#5}%
  \selectfont}%
\fi\endgroup%
\begin{picture}(1460,432)(1973,-1464)
\end{picture}%
 (resp. under the line 
\begin{picture}(0,0)%
\includegraphics{amoins.pstex}%
\end{picture}%
\setlength{\unitlength}{2072sp}%
\begingroup\makeatletter\ifx\SetFigFont\undefined%
\gdef\SetFigFont#1#2#3#4#5{%
  \reset@font\fontsize{#1}{#2pt}%
  \fontfamily{#3}\fontseries{#4}\fontshape{#5}%
  \selectfont}%
\fi\endgroup%
\begin{picture}(1460,432)(1973,-980)
\end{picture}%
) 
and $A_{+,-}(\sigma)$
is the number of pairs of edges such that one is above the line and
the other under and such that their supports do not intersect \begin{picture}(0,0)%
\includegraphics{aplusmoins.pstex}%
\end{picture}%
\setlength{\unitlength}{2072sp}%
\begingroup\makeatletter\ifx\SetFigFont\undefined%
\gdef\SetFigFont#1#2#3#4#5{%
  \reset@font\fontsize{#1}{#2pt}%
  \fontfamily{#3}\fontseries{#4}\fontshape{#5}%
  \selectfont}%
\fi\endgroup%
\begin{picture}(3170,474)(1973,-838)
\end{picture}%
.

For  $\sigma=(4,7,3,6,2,1,5)$, 
the pairs of edges
contributing to  $A_{+}(\sigma)$ are $\{(1,4),(3,3)\}$, $\{(2,7),(3,3)\}$ and $\{(2,7),(4,6)\}$.
The pair of edges
contributing to  $A_{-}(\sigma)$ is $\{(5,2),(6,1)\}$.
The pairs of edges
contributing to  $A_{+,-}(\sigma)$ are $\{(7,5),(1,4)\}$ and $\{(7,5),(3,3)\}$.

\begin{definition}
The number of nestings of a permutation $\sigma$ is equal to
$$
A_+(\sigma)+A_-(\sigma).
$$
\end{definition}

\begin{definition}
The number of alignments of a permutation $\sigma$ is equal to
$$
A_+(\sigma)+A_-(\sigma)+A_{+,-}(\sigma).
$$
\end{definition}
This definition looks a bit different but is equivalent to
the definition in \cite{wi}.

A descent in a permutation is an index $i$ such that $\sigma(i)>\sigma(i+1)$.
An ascent in a permutation is an index $i$ such that $\sigma(i)<\sigma(i+1)$.
We denote by $DES(\sigma)$ the number of descents of $\sigma$.

The pattern $31\mathrm{-}2$ (resp. $2\mathrm{-}31$, $13\mathrm{-}2$) occurs in $\sigma$ if there exist $i<j$ such that $\sigma(i)>\sigma(j)>\sigma(i+1)$ (resp.  $\sigma(j+1)<\sigma(i)<\sigma(j)$, $\sigma(i+1)>\sigma(j)>\sigma(i)$).
We denote by $(31\mathrm{-}2)(\sigma)$ the number of occurrences of the pattern
$31\mathrm{-}2$. For example, $(31\mathrm{-}2)(4,7,3,6,2,1,5)=2$, as
$73\mathrm{-}6$ and $73\mathrm{-}5$ are the two occurrences of the pattern $31\mathrm{-}2$.
\\

We will use bijections between permutations and weighted
bicolored Motzkin paths.
A bicolored Motzkin path of length $n$ is a sequence
$c=(c_1,\ldots ,c_n)$ such that 
$c_i\in\{N,S,E,\bar{E}\}$ for $1\le i\le n$ and such that
if 
$
h_i=\{ j<i\ |\ c_j=N\}-\{ j<i\ |\ c_j=S\}
$
then $h_1=0$, $h_i\ge 0$ for $2\le i\le n$ and $h_{n+1}=0$.

\section{Crossings, nestings and alignments}

\subsection{Crossings and nestings}

We use a bijection of Foata and Zeilberger \cite{fz} between permutations and weighted bicolored
Motzkin paths. We could also use
the bijection of Biane \cite{bi}. We refer to \cite{csz} for a compact definition
of these bijections.

To any permutation $\sigma$ we associate a pair $(c,w)$ made of a bicolored Motzkin path
$c=(c_1,c_2,\ldots ,c_n)$ and a weight $w=(w_1,\ldots ,w_n)$.
The path is created using the following rules~:
\begin{itemize}
\item $c_i=N$ if $i<\sigma(i)$ and $i<\sigma^{-1}(i)$
\item $c_i=E$ if $i\le \sigma(i)$ and $i\ge \sigma^{-1}(i)$
\item $c_i=\bar{E}$ if $i>\sigma(i)$ and $i<\sigma^{-1}(i)$
\item $c_i=S$ if $i>\sigma(i)$ and $i>\sigma^{-1}(i)$
\end{itemize}

The weight is created using the following rules~:
\begin{itemize}
\item $w_i=y p^{|A_+(i)|}q^{|C_+(i)|}$ if $c_i=N,E$.
\item $w_i=p^{|A_-(i)|}q^{|C_-(i)|}$ if $c_i=S,\bar{E}$.
\end{itemize}
This implies that $\prod_{i=1}^n w_i$ is equal to
$y^{WEX(\sigma)}q^{C_+(\sigma)+C_-(\sigma)}p^{A_+(\sigma)+A_-(\sigma)}$.
Let ${\mathcal P}_n$ be the set of pairs $(c,w)$ obtained from permutations
of $[n]$.
Then the coefficient of $y^k q^\ell p^m$ in
$$
\sum_{(c,w)\in {\mathcal P}_n} \prod_{i=1}^n w_i
$$
is $B(n,k,\ell,m)$ the number of permutations of $[n]$ with $k$
weak exceedances, $\ell$ crossings and $m$ nestings.

For example $\sigma=(4,1,5,6,2,3)$ gives the path $(N,\bar{E},N,E,S,S)$
and the weight $(y,1,yq,yq^2,q,1)$. \\

We now compute the generating function 
$\sum_n x^n \sum_{(c,w)\in {\mathcal P}_n} \prod_{i=1}^n w_i.
$
\begin{lemma}
If $i\le \sigma(i)$ then
$$
|C_+(i)|=h_i-|A_+(i)|;
$$
and if $i>\sigma(i)$ then
$$
|C_-(i)|=h_i-1-|A_-(i)|.
$$
\end{lemma}
\noindent{\bf Proof.} It is easy to prove by induction that
$$
h_i=|\{j<i\ |\ \sigma(j)\ge i\}|=|\{j\ge i\ |\ \sigma(j)<i\}|.
$$
Using the definitions given in Section 1, 
if $i\le \sigma(i)$ then
$$
A_+(i)\cup C_+(i)=\{j<i\ |\ \sigma(j)\ge i\}
$$
and if $i>\sigma(i)$ then
$$
A_-(i)\cup C_-(i)=\{j>i\ |\ \sigma(j)<i\}=\{j\ge i\ |\ \sigma(j)<i\}\backslash\{i\}.
$$
This implies the result. \qed

Thanks to that Lemma, using the machinery developed in \cite{f,vi}, we get directly that 
the generating function
$$
\sum_n x^n \sum_{(c,w)\in {\mathcal P}_n} \prod_{i=1}^n w_i.
$$
is
$$
\cfrac{1}{1-b_0x-\cfrac{\lambda_1x^2}{1-b_1x-\cfrac{\lambda_2x^2}
{1-b_2x-\cfrac{\lambda_3x^2}{\ddots}}}}
$$
with $b_n=y[n+1]_{p,q}+[n]_{p,q}$ and $\lambda_n=y[n]_{p,q}^2$.
This is $F(q,p,x,y)$ and proves Proposition \ref{th1}.


Note that this bijection implies that
\begin{proposition}
The number of permutations with $k$ weak exceedances and
$\ell$ crossings and $m$ nestings is equal to the number of permutations with $k$ weak exceedances and
$\ell$ nestings and $m$ crossings.
\end{proposition}

Similar results are known for set partitions and matchings \cite{st,kz,kl}.

\subsection{Link with alignments}

From the preceding result, we know that
the coefficient of $x^ny^kq^\ell$ in  $F(q,1,x,y)$ is
the number of permutations of $[n]$ with $k$ weak 
exceedances and $\ell$ crossings.

From Section 1 we also know that 
the number of permutations of $[n]$ with $k$ weak exceedances
and $\ell$ alignments is the coefficient of
$q^{(k-1)(n\mathrm{-}k)-\ell} x^n y^k$ in $\hat{E}(q,x,y)$. 

We want to prove Theorem \ref{big} which states that $\hat{E}(q,x,y)=F(q,1,x,y)$.
It follows from the following proposition:
\begin{proposition}
For any permutation with $k$ weak exceedances
the number of crossings plus the number of alignments is
$(k-1)(n\mathrm{-}k).$
\label{main}
\end{proposition}

{\noindent \bf Proof.} We suppose that $\sigma$ is a permutation of $[n]$
with $k$ weak exceedances.
For any $i$ with $1\le i<n$, we first define~:
\begin{itemize}
\item $B_+(i)=\{j\ |\ j<i\le \sigma(j)\}$  
\item $B_-(i)=\{j\ |\ \sigma(j)< i\le j\}$  
\end{itemize}
Note that  $h_i=|B_+(i)|=|B_-(i)|$ and that for $i>\sigma(i)$, 
$
A_-(i)\cup C_-(i)\cup \{i\}=B_-(i).$ \\
Therefore
\begin{equation}
A_-(\sigma)+C_-(\sigma)=\sum_{i>\sigma(i)} (|B_-(i)|-1)=k-n+\sum_{i>\sigma(i)} |B_+(i)|
\label{eq1}
\end{equation}

For $i\le \sigma(i)$, let
$$
E_+(i)=\{j\ |\ i\in C_+(j)\}=\{j\ | \ i<j\le\sigma(i)<\sigma(j)\}.
$$ 
It is easy to see that~:
$$
E_+(i)\cup A_+(i)=B_+(\sigma(i)).
$$

Therefore
\begin{eqnarray*}
A_+(\sigma)+C_+(\sigma)&=& \sum_{i\le \sigma(i)}|{C}_+(i)|+ |A_+(i)| \\
&=& \sum_{i\le \sigma(i)}|E_+(i)|+ |A_+(i)| \\
&=&\sum_{i\le \sigma(i)}|B_+(\sigma(i))| \\
&=&\sum_{i>\sigma(i)} |D_+(i)| \label{eq2},
\end{eqnarray*}
where
$$
D_+(i)=\{j\ | \ j\le \sigma(j)\ {\rm and}\ i\in B_+(\sigma(j))\}.
$$
It is easy to see that for $i>\sigma(i)$, 
$D_+(i) =\{j\ |\ j\le \sigma(j)\ {\rm and}\ \sigma(i)<\sigma(j)<i\}$ 
and therefore that
$
B_+(i)\cup D_+(i)\cup A_{+,-}(i)=\{j\ |\ j\ge \sigma(j)\}.
$
As they are pairwise disjoint then
$$
|B_+(i)|+ |D_+(i)|+ |A_{+,-}(i)|=k. 
$$

Combining equations (\ref{eq1}) and (\ref{eq2})
we get that
$A_+(\sigma)+C_+(\sigma)+A_-(\sigma)+C_-(\sigma)+A_{+,-}(\sigma)=
k-n+ \sum_{i>\sigma(i)} |D_+(i)|+|B_+(i)|+|A_{+,-}(i)|.$
This concludes
the proof of Proposition \ref{main}. \qed

\section{Permutation patterns}

Continued fractions like the one presented in Equation (\ref{best})
of Theorem \ref{max} were studied combinatorially in \cite{cm,csz}.
Theorem 10 in \cite{csz} associated with Theorem 22 in \cite{cm}
tells us that
the coefficient of $x^n y^k q^\ell p^m$ in
$F(x,y,q,p)$
is the number $D(n,k,\ell,m)$ of permutations $\sigma$ of $[n]$ with $n\mathrm{-}k$ descents, 
$\ell$ occurrences of the patterns $31\mathrm{-}2$ 
and $m$ occurrences of the pattern $2\mathrm{-}31$. This is Proposition \ref{th2}.
 
This can also be proved bijectively thanks to a bijection of Fran\c{c}on
and Viennot \cite{fv}. We present now that bijection. See also \cite{csz}.

Given a permutation $\sigma=(\sigma(1),\ldots , \sigma(n))$, we set $\sigma(0)=0$ and 
$\sigma({n+1})=n+1$. Let $\sigma(j)=i$.
Then $i$ is
\begin{itemize}
\item a valley if $\sigma({j-1})>\sigma({j})<\sigma({j+1})$
\item a double ascent  if $\sigma({j-1})<\sigma({j})<\sigma({j+1})$
\item a double descent  if $\sigma({j-1})>\sigma({j})>\sigma({j+1})$
\item a peak  if $\sigma({j-1})<\sigma({j})>\sigma({j+1})$
\item the beginning (resp. end) of a descent if $\sigma({j})>\sigma({j+1})$ 
(resp. $\sigma({j-1})>\sigma({j})$)
\item the beginning (resp. end) of a ascent if $\sigma({j})<\sigma({j+1})$ 
(resp. $\sigma({j-1})<\sigma({j})$).
\end{itemize}
If $\sigma(j)=i$ then we define $(31\mathrm{-}2)(i)$ (resp. $(2\mathrm{-}31)(i)$) to be the number of indices $k<j$ (resp. $k>j$)
such that $\sigma(k-1)>\sigma(j)>\sigma(k)$.\\

To any permutation, we associate a pair $(c,w)$ made of a bicolored Motzkin path
$c=(c_1,c_2,\ldots ,c_n)$ and a weight $w=(w_1,\ldots ,w_n)$.
The path is created using the following rules~:
\begin{itemize}
\item $c_i=N$ if $i$ is a valley
\item $c_i=E$ if $i$ is a double ascent
\item $c_i=\bar{E}$ if $i$ is a double descent
\item $c_i=S$ if $i$ is a peak
\end{itemize}

The weight is created using the following rules~:
\begin{itemize}
\item $w_i=y p^{(31\mathrm{-}2)(i)}q^{(2\mathrm{-}31)(i)}$ if $c_i=N,E$.
\item $w_i= p^{(31\mathrm{-}2)(i)}q^{(2\mathrm{-}31)(i)}$ if $c_i=S,\bar{E}$.
\end{itemize}
This implies that $\prod_{i=1}^n w_i$ is equal to
$y^{DES(\sigma)}p^{(31\mathrm{-}2)(\sigma)}q^{(2\mathrm{-}31)(\sigma)}$.

For example $\sigma=(6,2,1,5,3,4)$ gives the path $(N,\bar{E},N,E,S,S)$
and the weight $(y,1,yp,yp,y,1)$.\\

Now we prove the following Lemma
\begin{lemma}
For any $i$
$$
({31\mathrm{-}2})(i)+(2\mathrm{-}31)(i)=\left\{ \begin{array}{ll}
h_i & {\rm if}\ i \ {\rm is}\ {\rm the}\ {\rm beginning}\ {\rm of }\ {\rm an}\ {\rm ascent}\\ 
h_i-1 & {\rm if}\ i \ {\rm is}\ {\rm the}\ {\rm beginning}\ {\rm of }\ {\rm a}\  {\rm descent}\\ 
\end{array}\right.
$$
\end{lemma}
\noindent{\bf Proof.} 
We prove this lemma by induction. If $i$ is equal to 1 then 
$i$ is the beginning of an ascent and $(31\mathrm{-}2)(1)+(2\mathrm{-}31)(1)=0=h_1$.
 If $i>1$ then  $(31\mathrm{-}2)(i)+(2\mathrm{-}31)(i)=(31\mathrm{-}2)(i-1)+(2\mathrm{-}31)(i-1)+v$ where $v$ is zero,
 one or minus one.
It is easy to see that $v$ is one if $i-1$ is the end of a descent
and $i$ is the beginning of an ascent, $v$ is minus one if $i-1$ is the end of an ascent
and $i$ is the beginning of a descent and 0 otherwise. That gives exactly the lemma.  \qed\\

Let ${\mathcal P}_n$ be the set of pairs $(c,w)$ obtained from permutations
of $[n]$.
Using the machinery developed in \cite{f,vi}, we get directly that 
$\sum_n x^n \sum_{(c,w)\in {\mathcal P}_n} \prod_{i=1}^n w_i$ equals
$F(q,p,x,y)$ and proves Proposition \ref{th2}.\\

\section{Bijective proof of Theorem \ref{max}}

Combining the bijection of Fran\c{c}on and Viennot and the inverse of bijection
of Foata and Zeilberger \cite{fz}, we get Theorem \ref{max}.
We propose the direct mapping very similar to a bijection
proposed in \cite{csz}. Starting from 
a permutation $\sigma$ with $k$ descents, 
$\ell$ occurrences of the patterns $31\mathrm{-}2$ 
and $m$ occurrences of the pattern $2\mathrm{-}31$, we form a permutation
$\tau$ with $n\mathrm{-}k$ weak exceedances and $\ell$ crossings
and $m$ nestings.\\

We first form two two-rowed arrays $\tau_-$ and $\tau_+$.
We create the permutation $\tau$ from the tableaux of $\tau_-$ and $\tau_+$.
For $i$ from 1 to $n$, 
if $i$ is in the $j^{th}$ entry of the first row of $\tau_-$ (resp. $\tau_+$)
then $\tau(i)$ is the $j^{th}$ entry of the second row of $\tau_-$ (resp. $\tau_+$).
\\

The first row of $\tau_-$ contains all the entries of $\sigma$
that are the beginning of a descent. They are sorted in increasing order. 
The second row of $\tau_-$ contains all the entries of $\sigma$
that are the end of a descent. They are sorted such that
if $i$ is in first row then $\tau(i)<i$ and $(31\mathrm{-}2)(i)$ in $\sigma$ is equal to $A_-(i)$ in $\tau$.
One can easily check that this is always possible  in a unique way. 
For each $i$ in the first row starting from the smallest, $\tau(i)$ is the $((31\mathrm{-}2)(i)+1)^{th}$
smallest entry of the second row that is not yet chosen. Note that this 
implies that $(2\mathrm{-}31)(i)$ in $\sigma$ is equal to $C_-(i)$ in $\tau$.\\

The first line of $\tau_+$ contains all the entries of $\sigma$
that are the beginning of an ascent and that are sorted in increasing order. 
The second line of $\tau_+$ contains all the entries of $\sigma$
that are not the end of a descent. They are sorted such that
if $i$ is in first row then $\tau(i)\ge i$ and $(2\mathrm{-}31)(i)$ in $\sigma$ is equal to $C_+(i)$ in $\tau$.
One can again easily check that this is always possible in a unique way and that this 
implies that $(31\mathrm{-}2)(i)$ in $\sigma$ is equal to $A_+(i)$ in $\tau$.\\

For example if $\sigma=(5,1,7,4,3,6,8,2)$ then the 
$2\mathrm{-}31$ sequence is 
$$((31\mathrm{-}2)(1),\ldots ,(31\mathrm{-}2)(8))=(0,1,1,1,0,1,0,0)$$ 
$$((2\mathrm{-}31)(1),\ldots ,(2\mathrm{-}31)(8))=(0,0,1,1,2,1,1,0).$$
Then 
$$
\tau_-=\left(\begin{array}{l}
4,5,7,8\\
2,1,3,4 \end{array}\right)
$$
and 
$$
\tau_+=\left(\begin{array}{l}
1,2,3,6\\
8,5,6,7
 \end{array}\right)
$$
Then $\tau=(8,5,6,2,1,7,3,4)$.
One can check that $C_+(1)=0$, $C_+(2)=0$, $C_+(3)=1$,
$C_-(4)=1$, $C_-(5)=2$, $C_+(6)=1$, $C_-(7)=1$, $C_-(8)=0$.\\

If $\sigma$ is the image of $\tau$, it is easy to see that
\begin{lemma}
\begin{eqnarray*}
WEX(\tau)&=&n-DES(\sigma)\\
C_+(\tau)+C_-(\tau)&=&(2\mathrm{-}31)(\sigma)\\
A_+(\tau)+A_-(\tau)&=&(31\mathrm{-}2)(\sigma).
\end{eqnarray*}
\end{lemma}
Therefore if $\sigma$ has $k$ descents, $\ell$ occurrences of $(2\mathrm{-}31)$
and $m$ occurrences of $(31\mathrm{-}2)$ then $\sigma$ has
$n\mathrm{-}k$ weak exceedances, $\ell$ crossings and $m$ nestings.\\

This map is easily reversible.\\

We conclude this Section  by proving a similar result.
Given a permutation $\sigma=(\sigma(1),\ldots ,\sigma(n))$.
Let $\pi=(\sigma(n),\ldots ,\sigma(1))$. If $\sigma$ has $k-1$ ascents
(or $n\mathrm{-}k$ descents) and $\ell$ occurrences of the pattern $31\mathrm{-}2$ (resp. $2\mathrm{-}31$), then $\pi$
has $k-1$ descents and $\ell$ occurrences of the pattern $2-13$ (resp. $13\mathrm{-}2$).
Therefore
\begin{proposition}
There is a one-to-one correspondance between permutations
with $k$ weak exceedances and $\ell$ crossings and $m$ nestings 
and permutations with $k-1$ descents and $\ell$ occurrences
of the pattern $13\mathrm{-}2$ and $m$ occurrences of 
the pattern $2-13$.
\label{pattern}
\end{proposition}

\section{Generalization for decorated permutations}

We can also derive the generating function of
$$
A_{k,n}(q)=\sum_{i=0}^{k-1}{n\choose i}E_{k,n-i}.
$$
These were introduced in \cite{wi}. 

The following corollary is an easy consequence of Theorem
\ref{big}.
Let $A(q,x,y)=\sum_{n,k}A_{k,n}(q)x^n y^k$.
\begin{corollary}
$$
A(q,x,y)=
\cfrac{1}{1-b_0x-\cfrac{\lambda_1x^2}{1-b_1x-\cfrac{\lambda_2x^2}
{1-b_2x-\cfrac{\lambda_3x^2}{\ddots}}}}
$$
with $b_n=(1+y)[n+1]_q$ and $\lambda_n=yq[n]_q^2$.
\end{corollary}

Lauren Williams  \cite{wi0} proved that 
$$
A(q,x,y)=\frac{-y}{1-q}+\sum_{i\ge 1}
\frac{y^i(q^{2i+1}-y)}{q^{i^2+i+1}(q^i-q^i[i+1]_qx+[i]_qxy)}
$$

It would be interesting to have a direct proof of the identity
of the formal series form and the continued fraction form
of these generating functions \cite{wi0}. \\

We can also interpret these results combinatorially.
In \cite{wi} the coefficient of $q^{(n\mathrm{-}k)k-\ell}$ in $A_{k,n}(q)$
is interpreted in terms of decorated permutations with $k$
weak exceedances and $\ell$ alignments . Let us define these notions.
Decorated permutations
are permutations where the fixed points are bicolored \cite{po}.
We color these fixed points by colors $\{+,-\}$.
We say that $i\le_+ \sigma(i)$ if $i<\sigma(i)$ or $i=\sigma(i)$
and $i$ is colored with color $+$. We say that 
$i\ge_- \sigma(i)$ if $i>\sigma(i)$ or $i=\sigma(i)$
and $i$ is colored with color $-$. 

For a decorated permutation $\sigma$ and $i$, let
\begin{itemize}
\item $A_+(i)=\{j\ |\ j<i\le_+ \sigma(i)<\sigma(j)\}$
\item $A_-(i)=\{j\ | \ j>i\ge_- \sigma(i)>\sigma(j)\}$
\item 
$A_{+,-}(i)=\{j\ |\ i\ge_- \sigma(i)>\sigma(j)\ge j\}\cup \{j\ |\ 
\sigma(j)\ge j>i\ge_- \sigma(i)\}$
\item $C_+(i)=\{j\ | \ i<j\le\sigma(i)<\sigma(j)\}$
\item $C_-(i)=\{j\ | \ j>i>\sigma(j)>\sigma(i)\}$
\end{itemize}
As for permutations we define ${\rm A}_{+}(\sigma)=\sum_i {\rm A}_{+}(i)$.

With these notions, we can again define the number of alignments (resp. nestings, crossings)
of a decorated permutation $\sigma$ as $A_+(\sigma)+A_-(\sigma)+
A_{+,-}(\sigma)$ (resp. $A_+(\sigma)+A_-(\sigma)$, $C_+(\sigma)+C_-(\sigma)$). 
The number of weak exceedances (resp. non-exceedances) of a decorated permutation
is the cardinality of the set $\{i\ |\ i\ge_+ \sigma(i)\}$ (resp. $\{i\ |\ i< \sigma(i)\}$ ).

Let
$$
A(q,p,x,y)=\cfrac{1}{1-b_0x-\cfrac{\lambda_1x^2}{1-b_1x-\cfrac{\lambda_2x^2}
{1-b_2x-\cfrac{\lambda_3x^2}{\ddots}}}}
$$
with $b_n=(1+y)[n+1]_{p,q}$ and $\lambda_n=yq[n]_{p,q}^2$.

A direct generalization of the bijection \`a la 
Foata-Zeilberger on decorated permutations gives~:
\begin{proposition}
The coefficient of $x^n y^k q^\ell p^m$ in
$A(q,p,x,y)$
is the number of decorated permutations $\sigma$ of $[n]$ with $k$ weak 
exceedances, $\ell$ is the sum of the number of crossings and the number of 
non-exceedances and $m$ nestings.
\end{proposition}

We can also make a direct link with the alignments~:
\begin{proposition}
For any decorated permutation $\sigma$ with $k$ weak exceedances
$$
A_+(\sigma)+A_-(\sigma)+C_+(\sigma)+C_-(\sigma)+A_{+,-}(\sigma)
+|\{j\ |\ j>\sigma(j)\}|=(n\mathrm{-}k)k.
$$
\end{proposition}
\noindent{\bf Proof.} The proof is omitted as it follows exactly
the same steps as the proof of Proposition \ref{main}. \qed \\

We could also derive a  bijection \`a la 
Fran\c{c}on-Viennot on decorated permutations to interpret the
$A_{k,n}(q)$ in terms of descents and permutations 
patterns. One way would be to define new decorated permutations
where $i$ is bicolored if and only if $i$ is a double ascent and
$(2\mathrm{-}31)(i)=0$.

\section{Link with the  PASEP model}

The  PASEP model \cite{der} 
consists of black particles entering a row of $n$
cells, each of which is occupied by a black particle or vacant. A
particle may enter the system from the left hand side, hop to the
right or to the left and leave the system from the right hand side, with the
constraint that a cell contains at most one particle. 
We will say that the empty cells are filled with white
particles $\whi$. A basic configuration is a row of $n$ cells, each
containing either a black $\bla$ or a white $\whi$ particle. 
Let $\basic_n$ be the
set of basic configurations of $n$ particles. We write these
configurations as though they are words of length $n$ in the language
$\{\whi,\bla\}^*$.

The  PASEP defines a Markov chain $P$ defined on $\basic_n$ with the
transition probabilities $\alpha$, $\beta$, 
and $q$.  The probability $P_{X,Y}$, of finding the system in state
$Y$ at time $t+1$ given that the system is in state $X$ at time $t$ is
defined by:
\begin{subequations}
\begin{itemize}
\item If $X=A\bla\whi B$ and $Y=A\whi\bla B$ then
\begin{equation}
P_{X,Y}=1/(n+1);\ \ \ \
P_{Y,X}=q/(n+1)
\end{equation}
\item If $X=\whi B$ and $Y=\bla B$ then
\begin{equation}
P_{X,Y}=\alpha/(n+1).
\end{equation}
\item If $X=B\bla$ and $Y= B\whi $ then
\begin{equation}
P_{X,Y}=\beta/(n+1).
\end{equation}
\item Otherwise $P_{X,Y}=0$ for $Y\neq X$ and $P_{X,X}=1-\sum_{X\neq
    Y}P_{X,Y}$.
\end{itemize}
\end{subequations}

\begin{figure}[ht!]
  \begin{center}
    \includegraphics[scale=0.5]{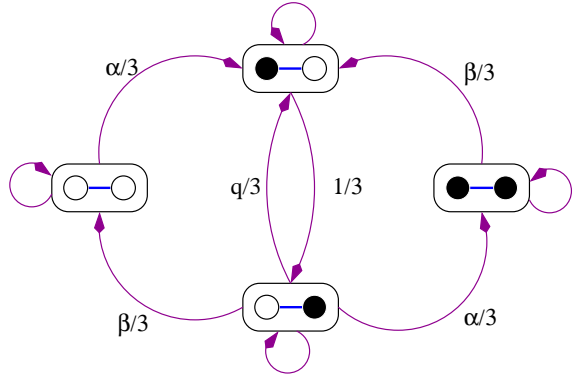}
  \end{center}\centering
  \caption{The chain $P$ for $n=2$.}
  \label{chainp}
\end{figure}

See an example for $n=2$ in Figure~\ref{chainp}.

This Markov chain has a unique stationary distribution \cite{der}.
The case $q=0$ was studied combinatorially by Duchi and Schaeffer
\cite{ds}. In \cite{bcer} another combinatorial approach was taken
to treat the general case.

\begin{definition}\cite{bcer}
  Let $\mathcal{P}(n)$ be the set of bicolored Motzkin paths of
  length $n$. The weight of the path in $\mathcal{P}(0)$ is 1.  The
  weight of any path $p$ denoted by $w(p)$ is the product of the weights of
  its steps. The weight of a step, $p_i$, starting at height $h$ is
  given by:
  \begin{subequations}
  \begin{eqnarray*}
    \mbox{if } p_i = N & \mbox{ then } & w(p_i) = [h+1]_q\\
    \mbox{if } p_i = \bar{E} & \mbox{ then } & w(p_i) =[h]_q+q^h/\alpha \\
    \mbox{if } p_i = E & \mbox{ then } & w(p_i) =[h]_q+q^h/\beta \\
    \mbox{if } p_i = S & \mbox{ then } & w(p_i) =[h]_q
+q^{h}/(\alpha\beta)-q^{h-1}(1/\alpha-1)(1/\beta-1).
  \end{eqnarray*}
\end{subequations}
\end{definition}

Given a path $p$, $\theta(p)$ is the basic configuration
such that each $\bar{E}$ and $S$ step is changed to $\whi$ and each
$E$ and $N$ step is changed to $\bla$. Let
\begin{equation}
  W(X) = \sum_{p \in \theta^{-1}(X)} w(p)
\end{equation}
and 
\begin{equation}
  Z_n = \sum_{X\in B_n} W(X)
\end{equation}

\begin{theorem}\cite{bcer}
At the steady state, the probability that the chain is in the basic 
configuration $X$ is
$$
\frac{W(X)}{Z_n}
$$
\end{theorem}

We can use these results and observations from the previous
sections to get~:
\begin{theorem}
If $\alpha=\beta=1$, at the steady state, the probability that 
the chain is in a basic configuration with $k$ particles is
$$
\frac{\hat{E}_{k+1,n+1}(q)}{Z_n}.
$$
\end{theorem}

Before proving that theorem, we need a Lemma
\begin{lemma}
There is a weight preserving bijection between $\mathcal{P}({n,k})$ the set of  
weighted bicolored Motzkin paths of length $n$ where the  weight  of any step starting at height $h$ is $[h+1]_q$
and where $k$ is the number of steps $N$ plus the number of steps
$E$ 
and  $\mathcal{P}'({n+1,k+1})$ the set of 
weighted bicolored Motzkin paths of length $n+1$
where the weight of any step starting at height $h$ is
 $[h+1]_q$ if the step is $N$ or $E$ and  $[h]_q$ otherwise
 and where $k+1$ is the number of steps $N$ plus the number of steps
$E$.  
\end{lemma}
\noindent{\bf Proof.} 

Starting with a weighted path $(c,w)$
in $\mathcal{P}({n,k})$ with $c=(c_1,\ldots ,c_n)$ and $w=(w_1,\ldots ,w_n)$, we create  a path
$(c',w')$ with
$c'=(c'_1,\ldots ,c'_{n+1})$ and $w'=(w'_1,\ldots ,w'_{n+1})$. 

We first set $c_0=N$ and $c_{n+1}=S$
and we  construct the path $c'$.
For $1\le i \le n+1$
\begin{itemize}
\item $c'_i=N$ if and only if $c_{i-1}=N$ or ${E}$ and $c_i=N$ or $\bar{E}$
\item $c'_i=E$ if and only if $c_{i-1}=N$ or ${E}$ and $c_i=S$ or ${E}$
\item $c'_i=\bar{E}$ if and only if $c_{i-1}=S$ or $\bar{E}$ and $c_i=N$ or $\bar{E}$
\item $c'_i=S$ if and only if $c_{i-1}=S$ or $\bar{E}$ and $c_i=S$ or ${E}$
\end{itemize}
It is easy to see that the path is a bicolored Motzkin path that it is of length $n+1$ 
and that $k+1$ is the number of steps $N$ plus the number of steps
$E$.  Moreover for $1\le i \le n$, if the starting height of $c_i$ is $j$ and $c'_i=N$ or $E$ 
(resp. $c'_i=S$ or $\bar{E}$) then
the starting height of $c'_i$ is $j$ (resp. $j+1$).
Therefore if we set $w'_i=w_i$ for $1\le i\le n$ and $w'_{n+1}=1$, then $(c',w')\in \mathcal{P}'({n+1,k+1})$.
This map is easily reversible. Remark : In this map the paths in $\mathcal{P'}$
do not have any steps $\bar{E}$ at height 0 and in particular
do not start or end with $\bar{E}$. This because
a path in $\mathcal{P'}$ having  a
$\bar(E)$ step at height 0 has weight zero and hence can be ignored.
\qed\\

\noindent{\bf Proof of the Theorem.} 
We want to compute for $\alpha=\beta=1$
$$
W({k,n})=\sum_{\begin{subarray}{c}X\in B_{n}\\
X \ {\rm has} \ k\ {\rm particles}\end{subarray}}W(X)=\sum_{p\in \mathcal{P}(n,k)}w(p).
$$
Now we use the previous lemma and get
$$
W({k,n})=\sum_{p\in \mathcal{P}'(n+1,k+1)}w(p).
$$
Using \cite{f,vi} and the definition of the weight of the steps of $\mathcal{P}'(n+1,k+1)$,
we  conclude that~:
$$
1+\sum_{n\ge 0} x^{n+1} \sum_{k=0}^{n} y^{k+1} W({k,n})=\hat{E}(q,x,y)
$$
and therefore that  $W({k,n})=\hat{E}_{k+1,n+1}(q)$.
\qed\\

\section{Conclusion}

Several open problems naturally arise from this work~:
\begin{itemize}
\item Can we define patterns for decorated permutations?
\item Can we generalize these $q$-Eulerian numbers to understand the  PASEP
when $\alpha\neq 1$ or $\beta\neq 1$? Can we use the permutations tableaux
\cite{sw,wi0,wi1}?
\item Can we extend the definition of $k$-crossing and $k$-nesting that were
defined for matchings and set partitions \cite{st} ?
\end{itemize}

\noindent{\bf Acknowledgments.} The author wants to thank Lauren Williams and
Petter Br\"and\'en for their encouragement and constructive comments and
the Institute Mittag Leffler where this work was done. The author also thanks the anonymous
referee for many useful comments.

\end{document}